\documentclass{article}
\usepackage{amsmath,amscd,amssymb}

\newtheorem{theorem}{Theorem}
\newtheorem{proposition}{Proposition}
\newtheorem{remark}{Remark}
\newtheorem{example}{Example}
\newtheorem{definition}{Definition}
\newtheorem{corollary}{Corollary}
\newtheorem{lemma}{Lemma}

\newcommand{\qed}{\hfill \hbox{\rule[-2pt]{3pt}{6pt}} \par}
\def\pd#1{ \partial_{#1} }

\def\nquad#1{\count0=0%
          \loop \ifnum#1 > \count0 \quad \advance \count0 by 1 \repeat}

\def\hsymbl#1{\smash{\hbox{\huge$#1$}}}
\def\hsymbu#1{\smash{\lower1.7ex\hbox{\huge$#1$}}}
\def\syz{\hbox{{\rm syz}}}

\title{\bf Incomplete Hypergeometric Systems Associated
to $1$-Simplex $\times$ $(n-1)$-Simplex}
\author{Kenta Nishiyama
 \footnote{Department of Mathematics, Kobe University and JST CREST.
}}
\date{December 17, 2010}

\begin{document}
\maketitle

\begin{abstract}
The ${\cal A}$-hypergeometric system was introduced by
Gel'fand, Kapranov and Zelevinsky in the 1980's.
Among several classes of ${\cal A}$-hypergeometric functions,
those for $1$-simplex $\times$ $(n-1)$-simplex are known to be a very nice class.
We will study an incomplete analog of this class.
\end{abstract}

\section{Introduction}
The ${\cal A}$-hypergeometric systems was introduced by
Gel'fand, Kapranov and Zelevinsky in the 1980's (\cite{GZK}).
It is a system of homogeneous differential equations with parameters
associated to an integer matrix $A$ and contains a broad class of hypergeometric
functions as solutions. 
Recently, the incomplete ${\cal A}$-hypergeometric system was proposed
toward applications to statistics and
a detailed study was given in the case of 
$ A=\left(
\begin{array}{cccc}
 1 & 1 & 0 & 0 \\
 0 & 0 & 1 & 1 \\
 0 & 1 & 0 & 1
\end{array}
\right) = $
$1$-simplex $\times$ $1$-simplex
(\cite{NT0}).
The system includes 
the incomplete Gauss' hypergeometric integral 
$\displaystyle I_{(a,b)}(\alpha,\beta,\gamma;x) = \int_a^b t^{\beta-1}(1-t)^{\gamma-\beta-1}(1-xt)^{\alpha}dt$
and the incomplete elliptic integral of the first kind
$\displaystyle F(z;k) = \int_0^z \frac{1}{\sqrt{(1-t^2)(1-k^2 t^2)}} dt$
as solution.
It is interesting to describe properties of these functions
in a general framework.
Among several classes of (complete) ${\cal A}$-hypergeometric functions,
those for $\Delta_1 \times \Delta_{n-1}$ ($1$-simplex $\times$ $(n-1)$-simplex)
are known to be a very nice class (see, e.g., \cite[Section 1.5]{SST}).

In this paper, we study an incomplete analog of this class.
In the section \ref{section:1n-hg}, we give a definition of an incomplete
$\Delta_1 \times \Delta_{n-1}$-hypergeometric system and 
prove that the existence of a solution of the system.
In the section \ref{section:series-sol},
we give a particular solution of the system and 
describe general solutions by combining with a base of the solutions
of (homogeneous) ${\cal A}$-hypergeometric system.
In the last section \ref{section:contiguity},
we give the complete list of contiguity relations for
the incomplete $\Delta_1 \times \Delta_{n-1}$-hypergeometric function.

\section{Incomplete $\Delta_1 \times \Delta_{n-1}$-hypergeometric system} \label{section:1n-hg}
We will work over the Weyl algebra in $2n$ variables
$D = {\bf C}\left \langle
\begin{array}{l}
 x_{11}, \cdots, x_{1n}, \pd{11}, \cdots, \pd{1n} \\
 x_{21}, \cdots, x_{2n}, \pd{21}, \cdots, \pd{2n}
\end{array}
\right \rangle $.

\begin{definition}\rm \label{def:hg-1n}
We call the following system of differential equations 
the {\it incomplete  $\Delta_1 \times \Delta_{n-1}$-hypergeometric system}:
\begin{equation} \label{eq:1n-system}
\left\{
\begin{array}{rll}
 (\theta_{i1}+\theta_{i2}-\alpha_i) \bullet f &= 0, &\quad (1 \leq i \leq n) \\
 \displaystyle \left (\sum_{i=1}^n \theta_{2i} + \gamma + 1 \right ) \bullet f &= [g(t,x)]_{t=a}^{t=b}, &\quad  \\
 (\pd{1i}\pd{2j}-\pd{1j}\pd{2i}) \bullet f &= 0, &\quad (1 \leq i < j \leq n)
\end{array}
\right.
\end{equation}
where $g(t,x) = t^{\gamma+1} \prod_{k=1}^{n} (x_{1k}+x_{2k}t)^{\alpha_k}$ and
$\alpha_i, \gamma \in {\bf C}$ are parameters.
The operator $\theta_{ij} = x_{ij} \pd{ij}$ is called the Euler operator.
\end{definition}

If $g(t,x)=0$ in (\ref{eq:1n-system}), the system agrees with
the ${\cal A}$-hypergeometric or GKZ hypergeometric system
associated to $\Delta_1 \times \Delta_{n-1}$.

\begin{remark} \rm
The incomplete $\Delta_1 \times \Delta_{n-1}$-hypergeometric system
introduced in Definition \ref{def:hg-1n} is a special but interesting case of
the incomplete ${\cal A}$-hypergeometric system (see appendix, \cite{NT0}).
Let $A$ be the following $(n+1) \times 2n$ matrix:
$$A 
=\left(
\begin{array}{ccccccc}
1 & 1 &   &   &        & \hsymbu{0}  &   \\
  &   & 1 & 1 &        &   &   \\
  &   &   &   & \ddots &   &   \\
  & \hsymbl{0}  &   &   &        & 1 & 1 \\
0 & 1 & 0 & 1 & \cdots & 0 & 1
\end{array}
\right).$$
We set 
$\beta = (\alpha_1, \ldots, \alpha_n, -\gamma-1 ) \in {\bf C}^{n+1}$
and $g = (0,\ldots,0,[g(t,x)]_{t=a}^{t=b})$.
Then the incomplete ${\cal A}$-hypergeometric system associated to $A, \beta, g$
is the incomplete $\Delta_1 \times \Delta_{n-1}$-hypergeometric system.
\end{remark}

We note that the ideal 
$\langle \pd{1i}\pd{2j}-\pd{1j}\pd{2i} \mid 1 \leq i < j \leq n \rangle$
generated by the third operators of (\ref{eq:1n-system})
is called the affine toric ideal associated to the matrix $A$
and it is denoted by $I_A$.
Moreover, $I_A$ is Cohen-Macaulay because $A$ is normal (\cite{Hochster}).

We note that the inhomogeneous system (\ref{eq:1n-system}) does not
necessarily have a solution $f$, when the inhomogeneous part
$[g(t,x)]_{t=a}^{t=b}$ is randomly given.

\begin{proposition}
For any $\alpha_i$, $\gamma \in {\bf C}$,
there exists a classical solution of the incomplete $\Delta_1 \times \Delta_{n-1}$-hypergeometric system.
\end{proposition}
{\it Proof.}
We may verify conditions (\ref{eq:syz1}) and (\ref{eq:syz2})
in Theorem \ref{theorem:syz} in the appendix with respect to $g=(0, \ldots, 0, [g(t,x)]_{t=a}^{t=b})$.
For $1 \leq i \leq n$, we have
\begin{align*}
 &\quad (\theta_{1i}+\theta_{2i}-\alpha_i) \bullet t^{\gamma+1} \prod_{k=1}^{n} (x_{1k}+x_{2k}t)^{\alpha_k} \\
 &= \alpha_i x_{1i} t^{\gamma+1} (x_{1i}+x_{2i}t)^{\alpha_i-1} \prod_{k \ne i}^{n} (x_{1k}+x_{2k}t)^{\alpha_k} \\
   &\qquad +\alpha_i x_{2i} t^{\gamma+2} (x_{1i}+x_{2i}t)^{\alpha_i-1} \prod_{k \ne i}^{n} (x_{1k}+x_{2k}t)^{\alpha_k}
   -\alpha_i t^{\gamma+1} \prod_{k=1}^{n} (x_{1k}+x_{2k}t)^{\alpha_k} \\
 &= \{ x_{1i} + x_{2i}t - (x_{1i}+x_{2i}t) \} \alpha_i t^{\gamma+1} (x_{1i}+x_{2k}t)^{\alpha_i-1} \prod_{k \ne i}^{n} (x_{1k}+x_{2k}t)^{\alpha_k} \\
 &= 0.
\end{align*}
Thus the condition (\ref{eq:syz1}) holds.

For $1 \leq i < j \leq n$, we have
\begin{align*}
 \pd{1i}\pd{2j} \bullet g(t,x)
 &= \pd{1i} \bullet \alpha_j t^{\gamma+2} (x_{1j}+x_{2j}t)^{\alpha_j-1} \prod_{k \ne j}^{n} (x_{1k}+x_{2k}t)^{\alpha_k} \\
 &= \alpha_i \alpha_j t^{\gamma+2} (x_{1i}+x_{2i}t)^{\alpha_i-1} (x_{1j}+x_{2j}t)^{\alpha_j-1} \prod_{k \ne i,j}^{n} (x_{1k}+x_{2k}t)^{\alpha_k}.
\end{align*}
Since this expression is symmetric in the indices $i$ and $j$,
we have $(\pd{1i}\pd{2j}-\pd{1j}\pd{2i}) \bullet g(t,x) = 0$.
Thus the condition (\ref{eq:syz2}) holds.
\qed
\bigbreak
Our definition of the incomplete $\Delta_1 \times \Delta_{n-1}$ hypergeometric system
is natural in terms of a definite integral with parameters.

\begin{proposition}
If ${\rm Re}\, \gamma$, ${\rm Re}\, \alpha_i > 0$, then the integral
\begin{equation}  \label{eq:intighg}
 \Phi(\beta; x) = \int_a^b t^\gamma \prod_{k=1}^{n} (x_{1k}+x_{2k}t)^{\alpha_k} dt
\end{equation}
is a solution of the incomplete $\Delta_1 \times \Delta_{n-1}$ hypergeometric system
(Definition \ref{def:hg-1n}).
\end{proposition} 

{\it Proof.}
From the general theory of ${\cal A}$-hypergeometric systems,
$\Phi(\beta; x)$ is annihilated by the elements of $I_A$ and
$\theta_{i1}+\theta_{i2}-\alpha_i$ for $1 \leq i \leq n$
(see, e.g., \cite[Section 5.4]{SST}).
We will prove that 
$$\left (\sum_{i=1}^{n} \theta_{2i} + \gamma + 1 \right ) \bullet \Phi(\beta; x)
= [g(t,x)]_{t=a}^{t=b}.$$
Applying $\sum_{i=1}^{n} \theta_{2i}$ to the integrand, we get
\begin{align*}
\left (\sum_{i=1}^{n} \theta_{2i} \right ) \bullet t^\gamma \prod_{k=1}^{n} (x_{1k}+x_{2k}t)^{\alpha_k}
 &= \sum_{i=1}^{n} \alpha_i x_{2i} (x_{1i}+x_{2i}t)^{\alpha_i-1} t^{\gamma+1} \prod_{k \ne i}^{n} (x_{1k}+x_{2k}t)^{\alpha_k} \\
 &= \sum_{i=1}^{n} t^{\gamma+1} \frac{\pd{} (x_{1i}+x_{2i}t)^{\alpha_i}} {\pd{} t} \prod_{k \ne i}^{n} (x_{1k}+x_{2k}t)^{\alpha_k} \\
 &= t^{\gamma+1} \frac{\pd{} \left (\prod_{k=1}^{n} (x_{1k}+x_{2k}t)^{\alpha_k} \right )} {\pd{} t}.
\end{align*}
By Stokes' theorem, we obtain
\begin{align*}
\left (\sum_{i=1}^{n} \theta_{2i} \right ) \bullet \Phi(\beta; x)
&= \int_a^b \left (\sum_{i=1}^{n} \theta_{2i} \right ) \bullet t^\gamma \prod_{k=1}^{n} (x_{1k}+x_{2k}t)^{\alpha_k} dt \\
&= \int_a^b t^{\gamma+1} \frac{\pd{} \left (\prod_{k=1}^{n} (x_{1k}+x_{2k}t)^{\alpha_k} \right )} {\pd{} t} dt \\
&= \left [ t^{\gamma+1} \prod_{k=1}^{n} (x_{1k}+x_{2k}t)^{\alpha_k} \right ]_{t=a}^{t=b} - (\gamma+1) \Phi(\beta; x).
\end{align*}
Thus the proposition is proved. \qed

\begin{example}\rm
We consider the following system of differential equations:
\begin{equation*} \label{eq:inhomogeneous11}
\left\{
\begin{array}{ll}
  (\pd{11} \pd{22} - \pd{12}\pd{21}) \bullet f & = 0, \\
  (\theta_{11} + \theta_{21} - \alpha_1) \bullet f & =  0, \\
  (\theta_{12} + \theta_{22} - \alpha_2) \bullet f & =  0, \\
  (\theta_{21} + \theta_{22} + \gamma + 1) \bullet f & = [g(t,x)]_{t=a}^{t=b}.
\end{array}
\right.
\end{equation*}
Here, $g(t,x) = t^{\gamma+1}(x_{11}+x_{21}t)^{\alpha_1}(x_{12}+x_{22}t)^{\alpha_2}$.

This is the incomplete $\Delta_1 \times \Delta_1$ hypergeometric system
for
$ A=\left(
\begin{array}{cccc}
 1 & 1 & 0 & 0 \\
 0 & 0 & 1 & 1 \\
 0 & 1 & 0 & 1
\end{array}
\right)$,
$\beta=(\alpha_1,\alpha_2, -\gamma-1)$, and
$g_1=0, g_2=0, g_3 = [g(t,x)]_{t=a}^{t=b}$.
\end{example}

A detailed study on the system is given in \cite{NT0}.

\section{Series Solution} \label{section:series-sol}
The Lauricella function $F_D$ is defined by
\begin{eqnarray*}
&& F_D(a,b_1, \ldots, b_n,c; z_1, \ldots, z_n) \\
&=& \sum_{m_1, \ldots, m_n = 0}^\infty 
         \frac{(a)_{m_1 + \cdots + m_n} (b_1)_{m_1} \cdots (b_n)_{m_n}}
              {(c)_{m_1 + \cdots + m_n} (1)_{m_1} \cdots (1)_{m_n}}
     z_1^{m_1} \cdots z_n^{m_n}.
\end{eqnarray*}
It is well-known that the Lauricella function $F_D$ of $n-1$ variables
gives a series solution of $\Delta_1 \times \Delta_{n-1}$-hypergeometric system.
We can give series solutions of our incomplete system 
in terms of the Lauricella series when parameters are generic.
We need $F_D$ of $n$ variables to give a solution.

\begin{theorem} \label{theorem:special-sol}
If $\gamma$ is not negative integer, the incomplete $\Delta_1 \times \Delta_{n-1}$-hypergeometric system
has a series solution which can be expressed in terms of 
the  Lauricella function $F_D$ as
\begin{align*}
F(\beta; x) = \prod_{k=1}^{n} x_{1k}^{\alpha_k} & \left ( \frac{b^{\gamma+1}}{\gamma+1}
F_D \left ( \gamma+1; -\alpha_1, \ldots, -\alpha_n; \gamma+2;
\frac{-x_{21}b}{x_{11}}, \ldots, \frac{-x_{2n}b}{x_{1n}} \right ) \right.\\
& \left.- \frac{a^{\gamma+1}}{\gamma+1}
F_D \left ( \gamma+1; -\alpha_1, \ldots, -\alpha_n; \gamma+2;
\frac{-x_{21}a}{x_{11}}, \ldots, \frac{-x_{2n}a}{x_{1n}} \right )
\right ).
\end{align*}
\end{theorem}

{\it Proof.}
For simplicity, we introduce some multi-index notations.
An $n$-dimensional multi-index is an $n$-tuple
$m=(m_1,\ldots,m_n)$
of non-negative integers.
The norm of a multi-index is defined by
$|m| = m_1 + \cdots + m_n$.
For a vector
$x_i = (x_{i1}, \ldots, x_{in})$ \,$(i=1,2)$,
define $x_i^m = x_{i1}^{m_1} \cdots x_{in}^{m_n}$
and for a vector
$\alpha=(\alpha_1,\ldots,\alpha_n) \in {\bf C}^n$,
define the Pochhammer symbol by
$(\alpha)_m = (\alpha_1)_{m_1} \cdots (\alpha_n)_{m_n}$.
By using these notations, the series $F$ can be written as
$$F = x_1^{\alpha} \sum_{m \geq 0}
c_m \left ( \frac{x_2}{x_1} \right )^m, \quad
c_m = \frac{(-1)^{|m|} (-\alpha)_m }{(\gamma+|m|+1) (1)_m}
(b^{\gamma+|m|+1}-a^{\gamma+|m|+1}). $$
We note that 
\begin{eqnarray*}
  \theta_{1k} \bullet F &=& (\alpha_k-m_k) F, \\
  \theta_{2k} \bullet F &=& m_k F.
\end{eqnarray*}
We now prove that the series $F$ satisfies the incomplete system (\ref{eq:1n-system}).
Firstly, $(\theta_{1i} + \theta_{2i} - \alpha_i) \bullet F = 0$ for $1 \leq i \leq n$
follows from above fact immediately.

Secondly, we will prove $\displaystyle \left (\sum_{i=1}^{n} \theta_{2i} + \gamma + 1 \right ) \bullet F = [g(t,x)]_{t=a}^{t=b}$,
which can be shown as
\begin{align*}
 \left (\sum_{i=1}^{n} \theta_{2i} + \gamma + 1 \right ) \bullet F
 &= (|m| + \gamma+1) F \\
 &= x_1^{\alpha} \sum_{m \geq 0}
    \frac{(-1)^{|m|} (-\alpha)_m }{(1)_m}
    (b^{\gamma+|m|+1}-a^{\gamma+|m|+1}) \left ( \frac{x_2}{x_1} \right )^m \\
 &= \left [ t^{\gamma+1} x_1^{\alpha} \sum_{m \geq 0}
    \frac{(-\alpha)_m }{(1)_m} \left ( -\frac{x_2 t}{x_1} \right )^m \right ]_{t=a}^{t=b} \\
 &= \left [ t^{\gamma+1} x_1^{\alpha} \prod_{k=1}^{n} \left ( 1+\frac{x_{2k} t}{x_{1k}} \right )^{\alpha_k} \right ]_{t=a}^{t=b} \\
 &= \left [ t^{\gamma+1} \prod_{k=1}^{n} \left ( x_{1k}+x_{2k} t \right )^{\alpha_k} \right ]_{t=a}^{t=b}.
\end{align*}
In the last two steps, we take a branch such that the equality holds.

Finally, we will prove $ (\pd{1i}\pd{2j}-\pd{1j}\pd{2i}) \bullet F = 0$
for $1 \leq i < j \leq n$.
This follows from the following two calculations:
\begin{align*}
\left (\theta_{1i}\theta_{2j}-\frac{x_{2j}x_{1i}}{x_{1j}{x_{2i}}}\theta_{1j}\theta_{2i} \right ) \bullet F
  &= x_1^{\alpha} \sum_{m \geq 0} (\alpha_i-m_i)m_j c_m \left ( \frac{x_2}{x_1} \right )^m \\
    & \qquad - x_1^{\alpha} \sum_{m \geq 0} (\alpha_j-m_j)m_i c_m \left ( \frac{x_2}{x_1} \right )^{m-e_i+e_j} \\
  &= x_1^{\alpha} \sum_{m \geq 0} (\alpha_i-m_i) (m_j+1) c_{m+e_j} \left ( \frac{x_2}{x_1} \right )^{m+e_j} \\
    & \qquad - x_1^{\alpha} \sum_{m \geq 0} (\alpha_j-m_j) (m_i+1) c_{m+e_i} \left ( \frac{x_2}{x_1} \right )^{m+e_j} 
\end{align*}
and
\begin{align*}
 (\alpha_i-m_i) (m_j+1) c_{m+e_j} 
 &= (\alpha_i-m_i) (m_j+1) \frac{(-1)^{|m+e_j|} (-\alpha)_{m+e_j}}{(\gamma+|m+e_j|+1) (1)_{m+e_j}} (b^{\gamma+|m+e_j|+1}-a^{\gamma+|m+e_j|+1})\\
 &= \frac{(-1)^{|m|+1} (-\alpha)_{m+e_j+e_i}}{(\gamma+|m|+2) (1)_{m}} (b^{\gamma+|m|+2}-a^{\gamma+|m|+2})\\
 &= (\alpha_j-m_j) (m_i+1) \frac{(-1)^{|m+e_i|} (-\alpha)_{m+e_i}}{(\gamma+|m+e_i|+1) (1)_{m+e_i}} (b^{\gamma+|m+e_i|+1}-a^{\gamma+|m_i|+1})\\
 &= (\alpha_j-m_j) (m_i+1) c_{m+e_i}.
\end{align*}
Therefore, the theorem is proved. \qed
\bigbreak
Gel'fand, Kapranov and Zelevinsky (\cite{GZK}) gave a base of 
the solutions of the (complete) ${\cal A}$-hypergeometric system.
We will give a base of solutions of our incomplete system by utilizing
their result and Theorem \ref{theorem:special-sol}.

For a parameter $\beta = (\alpha_1, \ldots, \alpha_n,-\gamma-1) \in {\bf C}^{n+1}$,
we set a $2 \times n$ matrix
$$
s^{(\ell)} = (s_{ij}^{(\ell)}) =
\left(
\begin{array}{ccccccc}
 \beta_1 & \cdots & \beta_{\ell-1} & \sum_{j=\ell}^{n+1} \beta_j & 0 & \cdots & 0 \\
 0 & \cdots & 0 & -\sum_{j=\ell+1}^{n+1} \beta_j & \beta_{\ell+1} & \cdots & \beta_{n}
\end{array}
\right)
$$
for $1 \leq \ell \leq n$.
Let $M^{(\ell)}$ be a set of $2 \times n$ matrices
$$M^{(\ell)} = \sum_{k=1}^{\ell-1} {\bf N}_0 \cdot (e_{2k}+e_{1\ell}-e_{1k}-e_{2\ell})
             + \sum_{k=\ell+1}^{n} {\bf N}_0 \cdot (e_{1k}+e_{2\ell}-e_{2k}-e_{1\ell}),$$
where $e_{ij}$ is the $2 \times n$ matrix 
whose $(i,j)$-entry is $1$ and the other entries are $0$.
We suppose that the condition of parameter $\beta$ called ``{\it $T$-nonresonant}",
that is the sets
$s^{(\ell)} \pm M^{(\ell)}$ $(1 \leq \ell \leq n)$
are pairwise disjoint (\cite[Definition 3]{GZK}).
Define series $\Psi^{(\ell)}(x)$ as
$$
 \Psi^{(\ell)}(x) = \Gamma(s^{(\ell)}+1) \sum_{k \in M^{(\ell)}} \frac{1}{\Gamma(s^{(\ell)}+k+1)} x^{s^{(\ell)}+k},
$$
where $\Gamma(s^{(\ell)}+k+1) = \prod_{i=1}^2 \prod_{j=1}^{n} \Gamma(s_{ij}^{(\ell)}+k_{ij}+1)$
and $x^{s^{(\ell)}+k} = \prod_{i=1}^2 \prod_{j=1}^{n} x_{ij}^{s_{ij}^{(\ell)}+k_{ij}}$.
These series are linearly independent and have the open domain
\begin{equation*} \label{eq:h-domain}
\left \vert \frac{x_{21}}{x_{11}} \right \vert < \left \vert \frac{x_{22}}{x_{12}} \right \vert < \cdots
< \left \vert \frac{x_{2n}}{x_{1n}} \right \vert
\end{equation*}
as a common domain of convergence.
Moreover they span the solution space of (complete)
$\Delta_1 \times \Delta_{n-1}$-hypergeometric system (\cite[Theorem 3]{GZK}).

By using this result, we obtain the following theorem.
\begin{theorem} \label{theorem:ser-sol}
Suppose the parameter $\beta$ is $T$-nonresonant and $\gamma$ is not negative integer.
\begin{enumerate}
\item
The common domain of convergence of $F(\beta;x)$ and $\Psi^{(\ell)}(x)$ is
$$ U : \left \vert \frac{x_{21}}{x_{11}} \right \vert < \left \vert \frac{x_{22}}{x_{12}} \right \vert < \cdots
< \left \vert \frac{x_{2n}}{x_{1n}} \right \vert < \frac{1}{\max(|a|,|b|)}.$$
\item
Any holomorphic solution of the incomplete system (\ref{eq:1n-system})
on $U$ can be written as
$$F(\beta;x) + \sum_{\ell=1}^{n} c_i \Psi^{(\ell)}(x), \qquad c_i \in {\bf C}.$$
\end{enumerate}
\end{theorem}

{\it Proof.}
Since the domain of convergence of the series $F$ is
$$ U_0 :
\left \vert \frac{x_{21}}{x_{11}} \right \vert < \frac{1}{\max(|a|,|b|)},
\left \vert \frac{x_{22}}{x_{12}} \right \vert < \frac{1}{\max(|a|,|b|)},
\cdots,
\left \vert \frac{x_{2n}}{x_{1n}} \right \vert < \frac{1}{\max(|a|,|b|)},
$$
we have the statement 1. The statement 2 is clear.
\qed
\bigbreak

\begin{theorem}
For $\sigma \in \mathfrak{S}_n$,
we suppose the parameter $\sigma(\beta)$ is $T$-nonresonant and $\gamma$ is not negative integer.
\begin{enumerate}
\item
The domain of convergence of the series $F(\beta;x)$ and $\sigma (\Psi^{(\ell)}(x))$ is
$$ \sigma(U) : \left \vert \frac{x_{2\sigma(1)}}{x_{1\sigma(1)}} \right \vert < \left \vert \frac{x_{2\sigma(2)}}{x_{1\sigma(2)}} \right \vert < \cdots
< \left \vert \frac{x_{2\sigma(n)}}{x_{1\sigma(n)}} \right \vert < \frac{1}{\max(|a|,|b|)}.$$
\item
Any holomorphic solution of the incomplete system (\ref{eq:1n-system})
on $\sigma(U)$ can be written as
$$F(\beta;x) + \sum_{\ell=1}^{n} c_i \sigma (\Psi^{(\ell)}(x)), \qquad c_i \in {\bf C}.$$
\end{enumerate}
Here, $\sigma (\Psi^{(\ell)}(x))$ is given
by the permutations
$x_{ij} \leftrightarrow x_{i \sigma(j)}$,
$s_{ij}^{(\ell)} \leftrightarrow s_{i \sigma(j)}^{(\ell)}$ and
$\beta_j \leftrightarrow \beta_{\sigma(j)}$.
\end{theorem}
{\it Proof.}
The theorem follows immediately from the $\sigma$-invariance of $F$.
\qed

\begin{remark} \rm
The closure of the union of $\sigma(U)$ coincides with the closure of $U_0$.
That is 
$$\overline{U_0} = \bigcup_{\sigma \in \mathfrak{S}_n} \overline{\sigma(U)}.$$
\end{remark}

\section{Contiguity Relation} \label{section:contiguity}
Contiguity relation is a relation among two functions
of which parameters are different by integer.
Miller (\cite{Miller}) gave contiguity relations for 
Lauricella functions and
Sasaki (\cite{Sasaki}) gave contiguity relations 
for Aomoto-Gel'fand hypergeometric functions, which include
the case of complete $\Delta_1 \times \Delta_{n-1}$-hypergeometric functions.
In \cite{SST2} and \cite{SST},
they give algorithms to compute contiguity relations
in the case of ${\cal A}$-hypergeometric functions.
These results are for complete functions.
In \cite{NT0},
an algorithm of computing contiguity relations under some conditions
is given for incomplete ${\cal A}$-hypergeometric functions and
the complete list of them for the incomplete
$\Delta_1 \times \Delta_1$-hypergeometric function is derived.

We will give contiguity relations of our incomplete system in this section.
We put $\delta = -\gamma -1$
(i.e., $\beta=(\alpha_1, \ldots, \alpha_n, \delta)$)
to make formulas of contiguity relations of 
the incomplete $\Delta_1 \times \Delta_{n-1}$-hypergeometric function simpler forms.
We put 
$$\Phi(\beta; x) = \int_a^b t^{-\delta-1} \prod_{k=1}^{n} (x_{1k}+x_{2k}t)^{\alpha_k} dt,$$
and assume ${\rm Re}\,(-\delta-1)$, ${\rm Re}\, \alpha_k > 0$.
Then, we note that $\Phi(\beta; x)$ is a solution of the following
incomplete $\Delta_1 \times \Delta_{n-1}$-hypergeometric system:
$$ 
\left\{
\begin{array}{cll}
 z_i \bullet f &= 0, & z_i := \theta_{1i}+\theta_{2i}-\alpha_i, \qquad (1 \leq i \leq n)\\
 z \bullet f &= [g(t,x)]_{t=a}^{t=b}, & z_{\,\,} := \sum_{i=1}^n \theta_{2i} - \delta, \\
 I_A \bullet f &= 0, &
\end{array}
\right.
$$
where 
$g(t,x) = t^{-\delta} \prod_{k=1}^{n} (x_{1k}+x_{2k}t)^{\alpha_k}$.
Let $a_{1k}$ and $a_{2k}$ be vectors corresponding to the $(2k-1)$-st and the $2k$-th columns of $A$ respectively.

\begin{theorem} \label{thm:contiguity}
The incomplete $\Delta_1 \times \Delta_{n-1}$-hypergeometric function
$\Phi(\beta ;x)$ satisfies the following contiguity relations.
\begin{itemize}
\item Shifts with respect to $a_{1k}$:
\begin{align}
S(\beta; -a_{1k}) \Phi(\beta; x) &= \alpha_{k} \Phi(\beta-a_{1k}; x), \label{eq:down-1k} \\ 
S(\beta-a_{1k}; +a_{1k}) \Phi(\beta-a_{1k}; x) &= \left (\sum_{i=1}^{n} \alpha_i - \delta \right ) \Phi(\beta; x) - [g(t,x)]_{t=a}^{t=b}, \label{eq:up-1k}
\end{align}
where
\begin{align*}
S(\beta; -a_{1k})  &= \pd{1k}, \\
S(\beta-a_{1k}; +a_{1k}) &= \sum_{i=1, i \ne k}^{n} (x_{1i}x_{2k}-x_{1k}x_{2i}) \pd{2i} + \sum_{i=1}^{n+1} \alpha_i x_{1k}.
\end{align*}

\item Shifts with respect to $a_{2k}$:
\begin{align}
S(\beta; -a_{2k}) \Phi(\beta; x) &= \alpha_{k} \Phi(\beta-a_{2k}; x), \label{eq:down-2k}\\
S(\beta-a_{2k}; +a_{2k}) \Phi(\beta-a_{2k}; x) &= \delta \Phi(\beta; x) + [g(t,x)]_{t=a}^{t=b}, \label{eq:up-2k}
\end{align}
where
\begin{align*}
S(\beta; -a_{2k})  &= \pd{2k}, \\
S(\beta-a_{2k}; +a_{2k}) &= \sum_{i=1, i \ne k}^{n} x_{1k}x_{2i}\pd{1i} + \left (\sum_{i=1, i \ne k}^{n} \theta_{2i} + \alpha_k \right ) x_{2k}.
\end{align*}
\end{itemize}
\end{theorem}
{\it Proof.}
The down-step relations (\ref{eq:down-1k}) and (\ref{eq:down-2k}) are easily verified.
We will prove only the up-step relations (\ref{eq:up-1k}) and (\ref{eq:up-2k}).

Let $L_1$ be the operator
$$\left ( \sum_{i=1, i \ne k}^{n} (x_{1i}x_{2k}-x_{1k}x_{2i}) \pd{2i} + \sum_{i=1}^{n} \alpha_i x_{1k} \right ) \pd{1k}
- \alpha_k \left ( \sum_{i=1}^{n} \alpha_i - \delta \right ) + \alpha_k \left ( \sum_{i=1}^{n} \theta_{2i} - \delta \right ).$$
We now prove that $L_1 \bullet \Phi(\beta; x) = 0$
which together with (\ref{eq:down-1k}) will prove the contiguity relation (\ref{eq:up-1k}).
The operator $L_1$ can be reduced by $z_k$ $(1 \leq k \leq n)$ as follows:
\begin{align*}
L_1 &= \sum_{i=1, i \ne k}^{n} x_{1i}x_{2k} \pd{2i} \pd{1k} - \sum_{i=1}^{n} (\theta_{2i} - \alpha_i) \theta_{1k} + \theta_{2k} \theta_{1k} 
          - \alpha_k \sum_{i=1}^{n} \alpha_i + \alpha_k \sum_{i=1}^{n} \theta_{2i} \\
  &= \sum_{i=1, i \ne k}^{n} x_{1i}x_{2k} \pd{2i} \pd{1k} - \sum_{i=1}^{n} (\theta_{2i} - \alpha_i) (\theta_{1k} + \theta_{2k} - \alpha_k) \\
  & \qquad \qquad \qquad + \sum_{i=1}^{n} (\theta_{2i} - \alpha_i) (\theta_{2k} - \alpha_k) + \theta_{2k} \theta_{1k}
          + \alpha_k \sum_{i=1}^{n} ( \theta_{2i} - \alpha_i ) \\
  &= \sum_{i=1, i \ne k}^{n} x_{1i}x_{2k} \pd{2i} \pd{1k} - \sum_{i=1}^{n} (\theta_{2i} - \alpha_i) z_k + \sum_{i=1}^{n} \theta_{2k} (\theta_{2i} - \alpha_i) + \theta_{2k} \theta_{1k} \\
  &= \sum_{i=1, i \ne k}^{n} x_{1i}x_{2k} \pd{2i} \pd{1k} - \sum_{i=1}^{n} (\theta_{2i} - \alpha_i) z_k + \sum_{i=1}^{n} \theta_{2k} (\theta_{1i} + \theta_{2i} - \alpha_i) - \sum_{i=1}^{n} \theta_{2k} \theta_{1i} + \theta_{2k} \theta_{1k} \\
  &= \sum_{i=1, i \ne k}^{n} x_{1i}x_{2k} (\pd{2i} \pd{1k} - \pd{2k} \pd{1i} ) - \sum_{i=1}^{n} (\theta_{2i} - \alpha_i) z_k + \sum_{i=1}^{n} \theta_{2k} z_i.
\end{align*}
Since the $\pd{2i} \pd{1k} - \pd{2k} \pd{1i}$ are elements of the toric ideal $I_A$,
we obtain $L_1 \bullet \Phi(\beta; x) = 0$.

Let $L_2$ be the operator
$$ \left (\sum_{i=1, i \ne k}^{n} x_{1k}x_{2i}\pd{1i} + \left (\sum_{i=1, i \ne k}^{n} \theta_{2i} + \alpha_k \right ) x_{2k} \right ) \pd{2k}
- \alpha_k \delta + \alpha_k \left ( \sum_{i=1}^{n} \theta_{2i} - \delta \right ).
$$
Since $L_2$ can be written as 
\begin{align*}
L_2 &= \sum_{i=1, i \ne k}^{n} x_{1k}x_{2i} (\pd{1i} \pd{2k} - \pd{2i} \pd{1k} ) + \sum_{i=1}^{n} \theta_{2i} z_k,
\end{align*}
we obtain $L_2 \bullet \Phi(\beta; x) = 0$
in an analogous calculation with the case of $L_1$.
\qed

Theorem \ref{thm:contiguity} gives contiguity relations
for $e_k = a_{1k} = (0,\ldots,0,\overset{k}{\check{1}},0,\ldots,0)$ $(1 \leq k \leq n)$,
but it does not give those for $e_{n+1} = (0,\ldots,0,1)$.
The set of vectors $\{e_1,\ldots,e_{n+1}\}$ is the standard basis of ${\bf Z}^{n+1}$.
The contiguity relations for $e_{n+1}$ can be obtained from Theorem \ref{thm:contiguity} as follows.

\begin{corollary}
The incomplete $\Delta_1 \times \Delta_{n-1}$-hypergeometric function
$\Phi(\beta ;x)$ satisfies the following contiguity relations.
\begin{itemize}
\item Shifts with respect to $e_{n+1}$:
\begin{align*}
S(\beta+e_{n+2}; -e_{n+2}) \Phi(\beta+e_{n+2}; x) &= \alpha_k \left (\sum_{i=1}^{n} \alpha_i - \delta \right ) \Phi(\beta; x) - \alpha_k [g(t,x)]_{t=a}^{t=b}, \\
S(\beta-e_{n+2}; +e_{n+2}) \Phi(\beta-e_{n+2}; x) &= \alpha_k \delta \Phi(\beta; x) + \alpha_k [g(t,x)]_{t=a}^{t=b},
\end{align*}
where
\begin{align*}
S(\beta+e_{n+2}; -e_{n+2})  &= S(\beta-a_{1k}; +a_{1k}) \pd{2k}, \\
S(\beta-e_{n+2}; +e_{n+2})  &= S(\beta-a_{2k}; +a_{2k}) \pd{1k} \text{ \quad for } 1 \leq k \leq n+1.
\end{align*}
\end{itemize}
\end{corollary}

Although we prove these contiguity relations for the integral representation
of the incomplete $\Delta_1 \times \Delta_{n-1}$ function,
they hold for functions which satisfy the system and
the two conditions (\ref{eq:down-1k}) and (\ref{eq:down-2k}).
By an easy calculation to check these conditions
for the series solution $F(\beta;x)$, we obtain the following corollary.

\begin{corollary}
The series solution $F(\beta;x)$ satisfies the same contiguity relations.
\end{corollary}

We note that $\Phi(\beta; x)$ can be formally expanded in $F(\beta;x)$.

\section{Appendix: A solvability of incomplete ${\cal A}$-hypergeometric systems} \label{section:appendix}
Let $D$ be the Weyl algebra in $n$ variables.
We denote by $A=(a_{ij})$ a $d\times n$-matrix whose elements are integers.
We suppose that the set of the column vectors of $A$ spans ${\bf Z}^d$.

\begin{definition} \rm  \label{def:ihg} (\cite{NT0})
We call the following system of differential equations $H_{A}(\beta,g)$
an {\it incomplete  ${\cal A}$-hypergeometric system}:
\begin{eqnarray*}
  (E_i-\beta_i)  \bullet f &=& g_i, \quad
  E_i-\beta_i= \sum_{j=1}^n a_{ij}  x_j \partial_{j} - \beta_i,
   \qquad(i = 1, \ldots, d)  \\
  \Box_{u,v} \bullet f &=& 0, 
\quad 
  \Box_{u,v} = \prod_{i=1}^n \partial_{i}^{u_{i}}
       - \prod_{j=1}^n \partial_{j}^{v_{j}}
\end{eqnarray*}
$$\mbox{ with } u, v \in {\bf N}_0^{n} \mbox{ running over all $u, v$ such that }  Au =  Av.$$
Here, ${\bf N}_0 = \{ 0, 1, 2, \ldots \}$, and 
$\beta = (\beta_1, \ldots, \beta_d) \in {\bf C}^d$ are parameters
and $g=(g_1, \ldots, g_d)$ where $g_i$ are given holonomic functions.
\end{definition}

We denote by $E-\beta$ the sequence $E_1-\beta_1, \ldots, E_d-\beta_d$
and $I_A$ the affine toric ideal generated by $\Box_{u,v}$ $(Au=Av)$
in ${\bf C}[\pd{1}, \ldots, \pd{n}]$.

\begin{lemma} \label{lemma:syz}
If the first homology of the Euler-Koszul complex vanishes;
$$H_1(K_{\bullet}(E-\beta; D/D I_A)) = 0,$$
then the syzygy module $\syz (E_1-\beta_1, \ldots, E_d-\beta_d) \subset (D/D I_A)^d$ 
is generated by $(E_i-\beta_i)e_j-(E_j-\beta_j)e_i$ $(1 \leq i < j \leq d)$.
\end{lemma}

{\it Proof.}
The Euler-Koszul complex of $D/D I_A$ is the following complex
$$
 0 \longrightarrow D/D I_A
   \stackrel{d_d}{\longrightarrow} \cdots
   \stackrel{d_3}{\longrightarrow} (D/D I_A)^{\binom{d}{2}}
   \stackrel{d_2}{\longrightarrow} (D/D I_A)^d 
   \stackrel{d_1}{\longrightarrow} D/D I_A
   \longrightarrow 0
$$
and the differential is defined by
$$ d_p(e_{i_1,\ldots,i_p}) 
= \sum_{k=1}^{p} (-1)^{k-1} (E_{i_k}-\beta_{i_k}) e_{i_1,\ldots,\widehat{i_k},\ldots,e_p}.$$
Here, $e_{i_1,\ldots,i_p}$ are basis vectors of $(D/D I_A)^{\binom{d}{p}}$.
The kernel of $d_1$ is $\syz(E_1-\beta_1, \ldots, E_d-\beta_d)$ 
and the image of $d_2$ is generated by $(E_i-\beta_i)e_j-(E_j-\beta_j)e_i$ $(1 \leq i < j \leq d)$
over $(D/D I_A)^d$.
Since the first homology is zero, the conclusion is obtained.
\qed

\begin{theorem}[\cite{Takayama}]  \label{theorem:syz}
If the first homology $H_1(K_{\bullet}(E-\beta; D/D I_A))$ vanishes
and the $g_i$ are holonomic functions satisfying
the following relations
\begin{eqnarray}
& &(E_i-\beta_i) \bullet g_j = (E_j-\beta_j) \bullet g_i, \qquad
   (i,j = 1, \ldots, d)  \label{eq:syz1} \\
& & \Box_{u,v} \bullet g_i = 0, \qquad
   (i=1, \ldots, d, Au=Av, u,v \in {\bf N}_0^n)  \label{eq:syz2}
\end{eqnarray}
then the incomplete hypergeometric system has a (classical) solution.
\end{theorem}

{\it Proof.}
By virtue of \cite[Theorem 4.1]{Kashiwara},
$\mathcal{E}xt_{\mathcal{D}}^1(\mathcal{D}/\mathcal{D} H_A(\beta),\mathcal{O})$ vanishes at generic points in ${\bf C}^n$.
Therefore, it is sufficient to prove that $\ell_1 g_1 + \cdots + \ell_d g_d = 0$
for all $(\ell_1,\ldots,\ell_d,\ell_{d+1},\cdots,\ell_{d+m}) \in \syz(E-\beta, \Box)$,
where $\Box$ is a finite sequence $\Box_{u_1,v_1}, \ldots, \Box_{u_m,v_m}$
which are generators of $I_A$.
Since for $(\ell_1,\ldots,\ell_d,\ell_{d+1},\cdots,\ell_{d+m}) \in \syz(E-\beta, \Box)$,
the relation $\sum_{i=1}^{d}\ell_i(E_i-\beta_i) + \sum_{i=1}^{m} \ell_{d+i}\Box_{u_i,v_i}=0$ holds,
we have $(\ell_1,\ldots,\ell_d) \in \syz(E-\beta)$ over $(D/D I_A)^d$.
By Lemma \ref{lemma:syz},
\begin{align*}
 \ell_1 g_1 + \cdots + \ell_d g_d
 &= (\ell_1, \ldots, \ell_d) \cdot g \\
 &= \sum_{1 \leq i < j \leq d} c_{ij} \{(E_i-\beta_i)e_j - (E_j-\beta_j)e_i\} \cdot g , \qquad c_{ij} \in {\bf C}\\
 &= \sum_{1 \leq i < j \leq d} c_{ij} \{(E_i-\beta_i)g_j - (E_j-\beta_j)g_i\} \\
 &= 0.
\end{align*}
\qed

\begin{remark} \rm
Matusevich, Miller and Walther (\cite[Theorem 6.3]{MMW})
showed that if the toric ideal $I_A$ is Cohen-Macaulay,
the $i$-th homology of the Euler-Koszul complex vanishes for all positive integers $i$.
\end{remark}
The following facts are known about Cohen-Macaulay property of toric ideals.
\begin{enumerate}
 \item If the initial monomial ideal of $I_A$ is square-free,
       then $A$ is normal (see, e.g., \cite[Proposition 13.15]{Sturmfels}).
 \item If the matrix $A$ is normal, then $I_A$ is Cohen-Macaulay (\cite{Hochster}).
\end{enumerate}
This is an easy tool for showing Cohen-Macaulayness of toric ideals.
When $A$ is $\Delta_1 \times \Delta_{n-1}$,
we can easily verify the condition 1.

\end{document}